\documentstyle[12pt]{article}
\def \reel{ {\rm I}\!{\rm R} }
 \newcommand{\too}{\longrightarrow}
\newcommand{\om}{\omega}
\newcommand{\Om}{\Omega}
\newcommand{\na}{\nabla}
\newcommand{\wi}{\widetilde}
\newcommand{\al}{\alpha}
\newcommand{\be}{\beta}
\newcommand{\ga}{\gamma}
\newcommand{\imp}{\Longrightarrow}

 \def \rat{ {\rm Q}\kern-.65em {}^{{}_/ }}

\newtheorem{Def}{Definition}[section]
\newtheorem{th}{Theorem}[section]
\newtheorem{pr}{Proposition}[section]
\newtheorem{Le}{Lemma}[section]
\title{Poisson manifolds with compatible pseudo-metric and pseudo-Riemannian
Lie algebras} \author{M. Boucetta} \date{}
\begin{document}
\maketitle

{\bf Abstract}\bigskip

The notion of Poisson manifold with compatible pseudo-metric was introduced
by the author in [1]. In this paper, we  introduce a new class of Lie
algebras which we  call a pseudo-Riemannian Lie algebras. The two notions are
strongly related: we prove that a linear Poisson structure on the dual of a
Lie algebra has a compatible pseudo-metric if and only if the Lie
algebra is a pseudo-Riemannian Lie algebra, and that the Lie algebra obtained by
linearizing at a point a Poisson manifold with compatible pseudo-metric 
 is a pseudo-Riemannian Lie
algebra. Furthermore, we  give some properties of the symplectic leaves of
such manifolds, and we prove that every 
Poisson manifold with compatible metric ( every Riemann-Lie
algebras) is unimodular. As a final, we  classify all pseudo-Riemannian Lie
algebras of dimension 2 and 3.

 \section{Introduction}
Let $P$ be a Poisson manifold with Poisson tensor $\pi$. A pseudo-metric
 of signature $(p,q)$ on the contangent bundle $T^*P$ is a smooth symmetric
contravariant 2-form $<,>$ on $P$ such that, at each point $x\in P$, $<,>_x$ is
non-degenerate on $T_x^*P$ with signature $(p,q)$. For any pseudo-metric $<,>$
on $T^*P$, we define a contravariant connection in the sens of Fernandes [2] and
Vaisman [4] by \begin{eqnarray*} 2<D_\al\be,\ga>&=&\#_{\pi}(\al).<\be
,\ga>+\#_{\pi}(\be).<\al ,\ga>-\#_{\pi}(\ga).< \al,\be>\\
&&+<[\al,\be]_{\pi},\ga >+<[\ga,\al]_{\pi}, \be>+<[\ga,\be]_{\pi}, \al>\qquad(1)
\end{eqnarray*}where $\al,\be,\ga\in\Om^1(P)$ and the Lie bracket $[,]_\pi$ is
given by \begin{eqnarray*}
[\al,\be]_{\pi}&=&L_{\#_{\pi}(\al)}\be-L_{\#_{\pi}(\be)}\al-d(\pi(\al,\be))\\
&=&i_{\#_{\pi}(\al)}d\be-i_{\#_{\pi}(\be)}d\al+d(\pi(\al,\be)).
\end{eqnarray*}
  $D$ will be called the Levi-Civita  contravariant connection associated with
the couple $(\pi,<,>)$. $D$ satisfies:

1. $D_\al\be-D_\be\al=[\al,\be]_{\pi};$

2. $\#_{\pi}(\al).<\be,\ga>=<D_\al\be,\ga>+<\be,D_\al\ga>.$

\begin{Def} With the notations above, the triple $(P,\pi,<,>)$ is called 
a pseudo-Riemannian Poisson manifold if, for any $\al,\be,\ga\in\Om^1(P)$,
$$D\pi(\al,\be,\ga)=
\#_{\pi}(\al).\pi(\be,\ga)-\pi(D_\al\be,\ga)-\pi(\be,D_\al\ga)=0.\eqno(2)$$When
$<,>$ is positive definite we call the triple a Riemann-Poisson manifold.
\end{Def}

Let $({\cal G},[,])$ be a Lie algebra and let $a$ be a non-degenerate bilinear
and symmetric form on $\cal G$. We define a bilinear map $A:{\cal G}\times{\cal
G}\too{\cal G}$ by
$$2a(A_uv,w)=a([u,v],w)+a([w,u],v)+a([w,v],u)\eqno(3)$$for any $u,v,w\in\cal G$.
$A$ satisifes:

1. $A_uv-A_vu=[u,v]$;

2. $a(A_uv,w)+a(u,A_uw)=0.$

\begin{Def} With the notations above, the triple $({\cal G},[,],a)$ is called
a pseudo-Riemannian Lie algebra if
$$[A_uv,w]+[u,A_wv]=0\eqno(4)$$for all $u,v,w\in\cal G$. When $a$ is positive
definite we call the triple $({\cal G},[,],a)$ a Riemann-Lie algebra.\end{Def}

In Section 1, we  give some basic properties of a pseudo-Riemannian Poisson
manifold and we  prove the following theorems:
\begin{th} Let $(P,\pi,<,>)$ be a pseudo-Riemannian Poisson manifold. Then, for
any point $x\in P$ such that the restriction of $<,>$ to $Ker\#_\pi(x)$ is
non-degenerate, the Lie algebra ${\cal G}_x$  obtained by linearizing the
Poisson structure at  $x$ is a  pseudo-Riemannian Lie algebra.\end{th}

\begin{th} Let $\cal G$ be a Lie algebra. The dual $({\cal G}^*,\pi)$ of $\cal G$
endowed with its linear Poisson structure $\pi$ has
a pseudo-metric $<,>$ for which the triple $({\cal G}^*,\pi,<,>))$ is a
pseudo-Riemannian Poisson manifold if and only if $\cal G$ is a
pseudo-Riemannian Lie algebra.\end{th}

{\bf Remark.} The condition on a  $x$ in  Theorem 1.1 is a condition
which depend only on the symplectic leaf of $x$ that means that if  $x$
satisfies this condition, since the parallel transport associated with $D$
preserve $Ker\#_\pi$ and $<,>$, each point on the symplectic leaf of $x$
satisfies the condition. Moreover, any point in which the Poisson tensor vanish
satisfies this condition and also any point in a Riemann-Poisson manifold.

In Section 2, we  prove the following theorems:
\begin{th}  Let $(P,\pi,<,>)$ be a Riemann-Poisson manifold and let $S$ be a
symplectic leaf of $P$. Then
$S$ is a K\"ahler manifold.
\end{th}
\begin{th} Let $(P,\pi,<,>)$ be a Riemann-Poisson manifold and let $S$ be a
symplectic regular leaf of $P$. Then the holonomy group of $S$ is finite.
\end{th}

In Section 3, we  prove the following theorem:
\begin{th} Every Riemann-Poisson manifold is unimodular. In particular, every
Riemann-Lie algebra is unimodular.\end{th}

  The  proof of the following theorem is a very long calculation and we don't
give it here.
 \begin{th} 1. The 2-dimensional abelian Lie algebra is the unique
2-dimensional pseudo-Riemann Lie algebra.

2. The 3-dimensional Lie algebras which has a pseudo-Riemannian Lie algebra
structure are:

a) The Heisenberg Lie algebra given by
$$[e_1,e_2]=e_3,\qquad [e_3,e_1]=[e_3,e_2]=0,$$

b) A family of Lie algebras given by
$$[e_1,e_2]=\al e_2+\be e_3,\quad [e_1,e_3]=\ga e_2-\al e_3,\quad [e_2,e_3]=0$$
where $\al,\be,\ga\in\reel$ and $\al^2+\be\ga\not=0$.

Furthermore, there is no Riemann-Lie algebra structure on the Heisenberg
Lie algebra, and a Lie algebra among the above familly has a structure
of Riemann-Lie algebra if and only if $\al^2+\be\ga<0$ and
$\ga>\be$.\end{th}

It seems, now,  that the condition of compatibility given by (3) is a nice
notion of compatibility beetwen Poisson structure and Riemann structure.
Moreover, the Riemann-Lie algebras, more than are examples of Riemann-Poisson
manifolds, give some new examples of homogeneous K\"ahler manifolds ( the
coadjoint orbits).

\section{Link beetwen pseudo-Riemannian Poisson manifolds and pseudo-Riemannian
Lie algebras}
\subsection{Basic materiel}
In this subsection, we collect the basic material which will be used throughout
this paper. Many fundamental definitions and results can be found in Vaisman's
monograph [4].

Let $P$ be a Poisson manifold with Poisson tensor $\pi$. The Poisson bracket on
$P$ is given by
$$\{f_1,f_2\}=\pi(df_1,df_2),\qquad f_1,f_2\in C^\infty(P).\eqno(5)$$
We also have a bundle map $\#_\pi:T^*P\too TP$ defined by
$$\be(\#_\pi(\al))=\pi(\al,\be),\qquad\al,\be\in T^*P.$$
On the space of differential 1-forms $\Om^1(P)$, the Poisson tensor induces a Lie
bracket
\begin{eqnarray*}
[\al,\be]_{\pi}&=&L_{\#_{\pi}(\al)}\be-L_{\#_{\pi}(\be)}\al-d(\pi(\al,\be))\\
&=&i_{\#_{\pi}(\al)}d\be-i_{\#_{\pi}(\be)}d\al+d(\pi(\al,\be)).\qquad\qquad\qquad\qquad\qquad
(6) \end{eqnarray*}

For this Lie bracket and the usual Lie bracket on vector fields, the bundle map
$\#_\pi$ induces a Lie algebra homomorphism $\#_\pi:\Om^1(P)\too{\cal X}(P)$:
$$\#_\pi([\al,\be]_{\pi})=[\#_\pi(\al),\#_\pi(\be)].\eqno(7)$$
We denote, as usual, by $X_f=\#_\pi(df)$ the hamiltonian associated with the
function $f\in C^\infty(P)$. We have
$$[df_1,df_2]_\pi=d\{f_1,f_2\}=L_{X_{f_1}}df_2,\qquad f_1,f_2\in
C^\infty(P).\eqno(8)$$

Let $<,>$ be a pseudo-metric on the contangent bundle $T^*P$ and $D$  the
Levi-Civita contravariant connection associated with the couple $(\pi,<,>)$.

For any $f,g,h\in C^\infty(P)$, we have
\begin{eqnarray*}
L_{X_f}<,>(dg,dh)&=&\#_\pi(df).<dg,dh>-<L_{X_f}dg,dh>-<dg,L_{X_f}dh>\\
&=&<D_{df }dg,dh>+<dg,D_{df }dh>\\
&&-<[df,dg]_\pi,dh>-<dg,[df,dh]_\pi>\\
&=&<D_{dg }df,dh>+<dg,D_{dh }df>.\end{eqnarray*}
So, we obtain, for any $f\in C^\infty(P)$ and any $\al,\be\in T^*P$,
$$L_{X_f}<,>(\al,\be)=<D_{\al }df,\be>+<\al,D_{\be }df>.\eqno(9)$$

Now, we  give a condition of compatibility beetwen the Poisson tensor $\pi$
and the pseudo-metric $<,>$ which is weaker than the condition given in the
Definition 1.1. This is possible because the Poisson tensor $\pi$ has vanishing
Schouten bracket $[\pi,\pi]_S$. In fact, since $D$ has vanishing torsion and
since the contravariant exterior differential $d_\pi$ associated with the bracket
$[,]_\pi$ is given by $d_\pi Q=-[\pi,Q]_S$, we can deduce obviously that, for any
$\al,\be,\ga\in\Om^1(P)$,
 $$-[\pi,\pi]_S(\al,\be,\ga)=D\pi(\al ,\be ,\ga
)+D\pi(\be ,\ga , \al)+ D\pi(\ga , \al, \be).\eqno(10)$$
From this formula and a straightforward calculation, we obtain, for any
$f,g,h\in C^\infty(P)$,
$$-[\pi,\pi]_S(df,dg,dh)-D\pi(dh, df, dg)=
\pi(D_{df}dh,dg)+\pi(df,D_{dg}dh).\eqno(11)$$

Now, we can give a condition of compatibility beetwen the Poisson tensor $\pi$
and the pseudo-metric $<,>$ which is weaker than the condition given by (2).
\begin{pr} Let $(P,\pi,<,>)$ be a Poisson manifold with a pseudo-metric on the
cotangent bundle and  $D$  the
Levi-Civita contravariant connection associated with the couple $(\pi,<,>)$.
The following assertions are equivalent:

1. The triple $(P,\pi,<,>)$ is a pseudo-Riemannian Poisson manifold.

2. For any $\al,\be\in\Omega^1(P)$ and any $f\in C^\infty(P)$,
$$\pi(D_\al df,\be)+\pi(\al,D_\be df)=0.\eqno(12)$$\end{pr} 

\subsection{ Proof of Theorem 1.1}

We recall the definition of the Lie algebra obtained by linearizing a Poisson
structure at a point. (See [6]).

Let  $(P,\pi)$ be a Poisson manifold and let   $x\in P$. We denote ${\cal
G}_x=Ker\#_\pi(x)$. For any $\al,\be\in{\cal G}_x$,  set
$$[\al,\be]_x=d_x(\pi(\wi\al,\wi\be))$$where $\wi\al,\wi\be\in\Omega^1(P)$
such that $\wi\al_x=\al$ and $\wi\be_x=\be.$
$({\cal G}_x,[,]_x)$ is a Lie algebra. 

Let $(P,\pi,<,>)$ be a pseudo-Riemannian Poisson manifold.
Fix a point $x\in P$ such that the restriction of $<,>$ to $Ker\#_\pi(x)$ is
non-degenerate. 
We denote $a$ the
restriction of  $<,>$ to ${\cal G}_x$ and $A$ the bilinear 2-form given  by
(3).

  For any $\al,\be\in{\cal G}_x$ and
$\wi\al,\wi\be\in\Omega^1(P)$ such that $\wi\al_x=\al$ and $\wi\be_x=\be$, we
claim that  $$(D_{\wi\al}\wi\be)_x=A_\al\be.$$
In fact, for each $\ga\in{\cal G}_x$ and
$\wi\ga\in\Omega^1(P)$ such that $\wi\ga_x=\ga$, we have 
\begin{eqnarray*}
2<D_{\wi\al}\wi\be,\wi\ga>&=&\#_{\pi}(\wi\al ).<\wi\be ,\wi\ga >+
\#_{\pi}(\wi\be ).<\wi\al ,\wi\ga >-\#_{\pi}(\wi\ga ).<\wi\al ,\wi\be >\\
&+&<[\wi\al ,\wi\be ]_\pi,\wi\ga>+<[\wi\ga ,\wi\al ]_\pi,\wi\be>+
<[\wi\ga ,\wi\be]_\pi,\wi\al>.\end{eqnarray*}
Since $\#_{\pi}(\wi\al )_x=\#_{\pi}(\wi\be )_x=\#_{\pi}(\wi\ga )_x=0$ and since 
$[\wi\al ,\wi\be ]_\pi(x)=[\al,\be]_x$, we have
$$2<(D_{\wi\al}\wi\be)_x,\ga>=2a(A_\al\be,\ga).$$

It remains to show that  $(D_{\wi\al}\wi\be)_x\in{\cal G}_x$. Indeed, for any 
$\mu\in\Omega^1(P)$, we have
\begin{eqnarray*}
\mu(\#_\pi(D_{\wi\al}\wi\be))&=&\pi(D_{\wi\al}\wi\be,\mu)\\
&=&\#_{\pi}(\wi\al ).\pi(\wi\be ,\mu)-\pi(\wi\be,D_{\wi\al}\mu)\\
&=&\#_{\pi}(\wi\al ).\pi(\wi\be ,\mu)-D_{\wi\al}\mu(\#_{\pi}(\wi\be )).
\end{eqnarray*}
Since $\#_{\pi}(\wi\al )_x=\#_{\pi}(\wi\be )_x=0$, we deduce
that  $(D_{\wi\al}\wi\be)_x\in{\cal G}_x$.

Now, let $f\in C^\infty(P)$ such that $d_xf=\ga.$ We have, from (12),
$$\pi(D_{\wi\al}df,\wi\be)+\pi(\wi\al,D_{\wi\be}df)=0.$$Differentiating this
relation at $x$, we obtain $$[A_\al\ga,\be]_x+[\al,A_\be\ga]_x=0$$and the theorem
follows.$\Box$

\subsection{Proof of Theorem 1.2}

Let $a$ be a bilinear non-degenerate symmetric form on $\cal G$ such that the
triple $({\cal G},[,],a)$ is a pseudo-Riemannian Lie algebra. We  define on
$T^*{\cal G}^*={\cal G}^*\times{\cal G}$ a pseudo-metric $<,>$ by setting
$$<(\mu,u),(\mu,v)>=a(u,v),\quad\mu\in{\cal G}^*,(u,v)\in{\cal G}^2.$$

For each vector $v\in\cal G$, we can define a linear form on ${\cal
G}^*$ denoted also $v$. For any $u,v\in\cal G$ and any $\mu\in{\cal
G}^*$,we have
$$\pi(dv,du)(\mu)=\mu([u,v]);\qquad [du,dv]_{\pi}=d[u,v];\qquad
D_{du}dv=d(A_uv). $$ 

It obvious that, in this case, (4) and (12) are equivalent.

Suppose that there is a pseudo-metric on ${\cal G}^*$ such that the
triple $({\cal G}^*,\pi,<,>)$ is a pseudo-Riemannian Poisson manifold. The
Lie algebra obtained by linearizing the Poisson structure at the origin of
${\cal G}^*$ is ${\cal G}$ and the theorem follows by Theorem 1.1.
$\Box$
\section{Sypmlectic leaves and the holonomy of Riemann-Poisson manifold}
This section is devoted to the proof of Theorem 1.3 and Theorem 1.4. Before, we
need some  Lemmas.
\begin{Le} Let $(P,\pi,<,>)$ be a pseudo-Riemannian Poisson manifold. Let
$S\subset P$ be a symplectic leaf an  $U$ an open of $P$. For any
$\al,\be\in\Om^1(P)$, we have
$$\#_{\pi}(\al)_{|S\cap
U}=0\qquad\mbox{or}\qquad\#_{\pi}(\be)_{|S\cap U}=0\Longrightarrow
\#_{\pi}(D_\al\be)_{|S\cap U}=0.\eqno(13)$$
$D$ is the Levi-Civita contravariant connection associated with $(\pi,<,>)$.
\end{Le}
{\bf Proof:} We have
$$\#_{\pi}(D_\al\be)-\#_{\pi}(D_\be\al)=[\#_{\pi}(\al),\#_{\pi}(\be)].$$
Hence $\#_{\pi}(\al)=0$ or $\#_{\pi}(\be)=0$ gives $
\#_{\pi}(D_\al\be)=\#_{\pi}(D_\be\al)$. 

Suppose that $\#_{\pi}(\be)=0$. For any $\ga\in\Om^1(P)$,
\begin{eqnarray*}
\ga(\#_{\pi}(D_\al\be))&=&\pi(D_\al\be,\ga)\\
&=&\#_{\pi}(\al).\pi(\be,\ga)-\pi(\be,D_\al\ga)=0
\end{eqnarray*}and the lemma follows.
$\Box$
\begin{Le} Let $(P,\pi,<,>)$ be a Riemann-Poisson manifold. Let
$O$ be the regular open on which the rank of the Poisson tensor is locally
constant. Then, we have:

1. $D$ is a $\cal F$-connection on $O$ in the sens of Fernandes [2], this means
that, for any $x\in O$ and any $\al\in T^*_xP$, we have
$$\#_\pi(\al)=0\qquad\imp\qquad D_\al=0;\eqno(14)$$

2. $D$ is a basic connection on $O$ in the sens of Fernandes [2], this means that
for any symplectic leaf $S\subset O$ and for any $\al,\be\in\Om^1(O)$ such that
$\#_\pi(\be)_{|S}=0$, we have
$$(D_\al\be)_{|S}={[\al,\be]_{\pi}}_{|S}.\eqno(15)$$\end{Le}
{\bf Proof:}  Let $U$ be an open in $O$ on which the rank of the Poisson tensor
is constant. On $U$, the cotangent bundle splits
$$T^*P=Ker\#_\pi\oplus Ker\#_\pi^{\perp}$$where $Ker\#_\pi^{\perp}$ is the
$<,>$-orthogonal of $Ker\#_\pi$.

Let $x\in U$ and  $\al\in T^*_xP$ such that $\#_\pi(\al)=0$. 
Choose a 1-forme $\wi\al$ on $U$ such that $\#_\pi(\wi\al)=0$ and
$\wi\al_x=\al$. For any $\be\in\Om^1(U)$ we have, by Lemma 3.1,
$\#_\pi(D_{\wi\al}\be)=0$. We claim that $D_{\wi\al}\be\in
Ker\#_\pi^{\perp}$. Indeed, for any $\ga\in Ker\#_\pi$, we have
$$<D_{\wi\al}\be,\ga>=\#_\pi(\wi\al).<\be,\ga>-<\be,D_{\wi\al}\ga>=
-<\be,D_{\wi\al}\ga>.$$Now, by the spliting of $T^*P$,
$\be=\be_1+\be_1^{\perp}$ and so $<\be,D_{\wi\al}\ga>=<\be_1,D_{\wi\al}\ga>$.
This quantity vanish  from the definition of $D$, so $D$ is a $\cal
F$-connexion. Since $D$ has a vanishing torsion, is also a basic
connection.$\Box$
\begin{Le}  Let $(P,\pi,<,>)$ be a Riemann-Poisson manifold. Let
$O$ be the regular open on which the rank of the Poisson tensor is locally
constant. Then, we have:

1. For any $x\in O$, for any $\al,\be\in Ker\#_\pi(x)$ and for any $f\in
C^\infty(O)$, we have 
$$L_{X_f}<,>(\al,\be)=0;$$

2. For any Casimir functions $f,g$, $<df,dg>$ is a Casimir function.
\end{Le}
{\bf Proof:} 1. follows from Lemma 3.2 and (9).

2. Let $(f,g)$ be a couple of Casimir functions on $P$ and let $h\in
C^\infty(P)$. We have
\begin{eqnarray*}
\{h,<df,dg>\}&=&X_h.<df,dg>\\
&=&<D_{dh}df,dg>+<df,D_{dh}dg>\\
&=&<D_{df}dh+[dh,df]_\pi,dg>+<df,D_{dg}dh+[dh,dg]_\pi>.\end{eqnarray*}
Now $[dh,df]_\pi=d\{h,f\}=0$ and $D_{df}dh$ vanish on $O$ by Lemma 3.2. So 
$\{h,<df,dg>\}$ is zero on $P$ since $O$ is dense.$\Box$
\subsection{Proof of Theorem 1.3}
Let $(P,\pi,<,>)$ be a Riemannian-Poisson manifold and let $S\subset P$ be a
symplectic leaf. We denote the symplectic form of $S$ by $\om_S$.

For any vectors fields $X,Y$ tangent to $S$, we set
$$\na_X^SY=\#_\pi(D_\al\be)_{|S}\eqno(16)$$where $\#_\pi(\al)_{|S}=X$ and
$\#_\pi(\be)_{|S}=Y$. It follows from Lemma 3.1 that $\na^S$ define a
torsionless covariant connection on $S$ and, obviously, we have
$$\na^S\om_S=0.$$

For any $x\in S$, we have
$$T^*_xP=Ker\#_\pi(x)\oplus(Ker\#_\pi(x))^{\perp}$$and so the linear map
$\#_\pi(x):(Ker\#_\pi(x))^{\perp}\too T_xS$ is an isomorphism.
For any $u,v\in T_xS$, we set
$$g_S(u,v)=<\#_\pi(x)^{-1}(u),\#_\pi(x)^{-1}(v)>.\eqno(17)$$$(S,g_s)$ is a
Riemann manifold and $\na^Sg_S=0$ and so $\na^S$ is the Levi-Civita connection
of $g_S$.

Now, it is classical (see [4]) that $TS$ has a $\na^S$-parallel complex structure
$J=A(-A^2)^{-\frac12}$, where $A$ is given by 
$$\om_S(u,v)=g_S(Au,v).$$ So $S$ is a K\"ahler manifold.$\Box$
\subsection{Proof of Theorem 1.4}
The linear Poisson holonomy was introduced by Ginzburg and Golubev in $[3]$. A
more deep study of this notion was done by Fenrnandes in $[2]$ essentially by the
notion of basic connection. In the same paper, Fernandes introduces the notion of
Poisson holonomy. For a regular leaf, the (linear)
Poisson holonomy coincides with  standart ( linear) holonomy.

Let $(P,\pi,<,>)$ be a Riemann-Poisson manifold and  $S$  a regular
symplectic leaf. We  prove that the holonomy group of $S$ is finite.

{\bf First step:} The linear holonomy group of $S$ is finite.

The linear holonomy group of $S$ coincides with the linear Poisson holonomy since
the leaf is regular. Now, the linear Poisson holonomy can be determined by a
basic connection ( see [R]). The Levi-Civita contravariant
connection $D$ is, by Lemma 3.2, a basic connection on a neighbourhood of $S$.
Moreover,  the parallel transport defined by  $D$ is given by isometries on
$Ker\#_\pi$. This gives the claim.

{\bf Second step:} The  holonomy group of $S$ is finite.

The holonomy group coincides with the  Poisson holonomy group since the
leaf is regular. The Poisson holonomy is given by  the hamiltonian flows. Now,
let $x\in S$ and $(q_1,\ldots,q_k,p_1,\ldots,p_k,y_1,\ldots,y_l)$ a	Darboux
coordinates defined on a neighbourhood $U$ of $x$. We considere the submanifold
$N$ of $U$   defined by $p=q=0$. We define a metric $g_N$ on $T^*N$
by $$g_N(dy_i,dy_j)=<dy_i,dy_j>.\eqno(17)$$ By considering $(N,g_N)$, we have
according to Lemma 3.3 a well-defined notion of "transverse Riemann structure"
along $S$ and the holonomy preserve this transverse structure. Now, the elements
of the holonomy group are isometries and their differentials at $x$ are elements
of linear holonomy group which is finite so the holonomy group is finite.$\Box$
\section{Poisson-Riemann manifolds are unimodular}
This section is devoted to the proof of Theorem 1.5.

We begin with a metric version of the modular vector field. For the details on
the modular vector field see $[5]$.

Let $(P,\pi)$ be a Poisson manifold and $g$ a Riemannian metric on $P$. We
denote $\#_g:T^*P\too TP$ the musical isomorphism associated with $g$. Also $g$
defines a density on $P$ which we denote $\mu_g$. The modular vector field with
respect to $\mu_g$ is given by
$$\phi_{\mu_g}(f)=div_g(X_f)=-\sum_{i=1}^{}g(\na_{e_i}X_f,e_i)\eqno(18)$$where
$\na$ is the Levi-Civita connection associated to $g$, $(e_1,\ldots,e_n)$ is a
local orthonormal basis of vector fields and $n=dimP$.

If $h(u,v)=g(Ju,v)$ is another Riemann metric, we have $\mu_h=\sqrt{detJ}\mu_g$
and
$$\phi_{\mu_h}=\phi_{\mu_g}-\frac12 X_{Ln(detJ)}.\eqno(19)$$

Now, we define a Riemann metric $<,>$ on the cotangent bundle $T^*P$ by
$$<\al,\be>=g(\#_g(\al),\#_g(\be)),\qquad\al,\be\in T^*P.$$

(We remark that any Riemann metric on the cotangent bundle $T^*P$
can be obtained in this may).

Let $D$ be the Levi-Civita contravariant connection associated with
$(\pi,<,>)$. We claim that
$$\phi_{\mu_g}=\frac12\sum_{i=1}^nL_{X_f}<,>(\al_i,\al_i)=
\sum_{i=1}^n<D_{\al_i}df,\al_i>\eqno(20)$$where $(\al_1,\ldots,\al_n)$ is a
local othonormal basis of 1-forms. This fact follows easily from the
following formula which is also easy to etablish:
$$L_{X_f}<,>(\al,\al)=-L_{X_f}g(\#_g(\al),\#_g(\be)),\qquad\al,\be\in
T^*P.\eqno(21)$$

Now, we considere a Riemann-Poisson manifold $(P,\pi,<,>)$ and denote
$\phi_{<,>}$ the modular vector field defined by $(20)$. We will show that 
$\phi_{<,>}$ is zero on $P$ by showing that it is zero on the regular open $O$.

Let $x\in O$ and $S$ the symplectic leaf of $x$. We considere the symplectic
form $\om_S$ of $S$, the Riemannian metric $g_S$ on $S$ given by $(18)$ and
$\na^S$ its Levi-Civita connection.

We have, in a neighbourhood of $x$, $$T^*P=Ker\#_\pi\oplus(Ker\#_\pi)^{\perp}.$$
We choose $(\al_1,\ldots,\al_l)$ and $(\be_1,\ldots,\be_{n-l})$ two local
orthonormal basis respectively of $Ker\#_\pi$ and $(Ker\#_\pi)^{\perp}$. We
have, according to Lemma 3.2,
\begin{eqnarray*}
\phi_{<,>}(f)(x)&=&\sum_{i=1}^{n-l}<D_{\be_i}df,\be_i>\\
&=&\sum_{i=1}^{n-l}g_S(\na^S_{\#_\pi(\be_i)}X_f,\#_\pi(\be_i))\\
&=&-div_{g_S}(X_f)(x).\end{eqnarray*}
Now, according to Theorem 1.3 there is a Riemann metric $h$ on $S$ such that
$(S,h,\om_S)$ is a K\"ahler manifold and such that the isomorphism $J$ given by
$h(u,v)=g_S(Ju,v)$ is $\na^S$-parallel. It follows that $detJ$ is constant and
$div_{g_S}(X_f)=div_{h}(X_f)$. To conclude, we recall the well-known fact that in
a K\"ahler manifold the divergence with respect to the K\"ahler metric of any
hamiltonian vector field is zero.

Now, let $\cal G$ be a Riemann-Lie algebra. According to Theorem 1.2, ${\cal
G}^*$ inherits a structure of Poisson-Riemann manifold which is unimodular and
the theorem follows.$\Box$\bigskip

{\bf References}

\bigskip

[1] M. Boucetta,  Compatibilit\'e des structures pseudo-riemanniennes et des
structures de Poisson, C. R. Acad. Sci. Paris, t. 333, S\'erie I, p. 763-768,
2001.

[2]  R. L. Fernandes, Connections in Poisson geometry 1: holonomy and
invariants, J. Diff. Geom. 54 , p. 303-366, 2000.

[3] V. Ginzburg and A. Golubev, Holonomy on Poisson Manifolds and The Modular
Class, Preprint math.DG/9812153.

[4] I. Vaisman, Lectures on the Geometry of Poisson Manifolds, Progress in
Mathematics, vol. 118, Birkh\"auser, Berlin, 1994.

[4] A. Weinstein, The modular Automorphism Group of a Poisson Manifold, J.
Geom. Phys. 23, p. 379-394, (1997).

[5] A. Weinstein, The Local Structure of Poisson Manifolds, J. Diffrential
Geometry 18, p. 523-557,1983.\bigskip

{\bf Mohamed BOUCETTA

Facult\'e des sciences et techniques, BP 549, Marrakech, Maroc

Email: boucetta@fstg-marrakech.ac.ma}

 \end{document}